\documentclass[12pt]{article}

\usepackage{amsmath, epsfig, cite}
\usepackage{amssymb}
\usepackage{amsfonts}
\usepackage{latexsym}
\usepackage{float}
\usepackage{color}

\newtheorem{thm}{Theorem}[section]

\newtheorem{conj}[thm]{Conjecture}
\newtheorem{lem}[thm]{Lemma}

\newcommand*{\pFq}[5]{{}_{#1}F_{#2}\left[ \begin{matrix} #3\\[5pt] #4\end{matrix};#5\right]}

\numberwithin{equation}{section}

\newcommand{\qed}{{\hfill$\square$}\medskip}

\begin{document}

\begin{center}
{\large\bf Proof of Sun's supercongruence involving Catalan numbers}
\end{center}

\vskip 2mm \centerline{Ji-Cai Liu}
\begin{center}
{\footnotesize Department of Mathematics, Wenzhou University, Wenzhou 325035, PR China\\
{\tt jcliu2016@gmail.com } \\[10pt]
}
\end{center}


\vskip 0.7cm \noindent{\bf Abstract.}
We confirm a conjectural supercongruence involving Catalan numbers, which is one of
the 100 selected open conjectures on congruences of Sun. The proof makes use of
hypergeometric series identities and symbolic summation method.


\vskip 2mm
\noindent{\it MR Subject Classifications}: 11A07, 11B65, 05A19, 33C20

\section{Introduction}
In 2003, Rodriguez-Villegas \cite{rv-b-2003} conjectured the following four supercongruences associated to certain elliptic curves:
\begin{align*}
\sum_{k=0}^{p-1}\frac{{2k\choose k}^2}{16^k}\equiv \left(\frac{-1}{p}\right) \pmod{p^2},\quad
\sum_{k=0}^{p-1}\frac{{2k\choose k}{3k\choose k}}{27^k}\equiv \left(\frac{-3}{p}\right) \pmod{p^2},\\
\sum_{k=0}^{p-1}\frac{{2k\choose k}{4k\choose 2k}}{64^k}\equiv \left(\frac{-2}{p}\right) \pmod{p^2},\quad
\sum_{k=0}^{p-1}\frac{{3k\choose k}{6k\choose 3k}}{432^k}\equiv \left(\frac{-1}{p}\right) \pmod{p^2},
\end{align*}
where $p\ge 5$ is a prime and $\left(\frac{\cdot}{p}\right)$ denotes the Legendre symbol.
These four supercongruences were first proved by Mortenson \cite{mortenson-jnt-2003,mortenson-tams-2003} by using the Gross-Koblitz formula.
Guo, Pan and Zhang \cite{gpz-jnt-2017} established some interesting $q$-analogues of 
the above four supercongruences. For more $q$-analogues of congruences, one can refer to \cite{gl-jdea-2018,gs-ca-2019,gs-rm-2019}.

Let $E_n$ be the $n$th Euler number.
In 2016, Z.-H. Sun \cite{sunzh-ijnt-2016} proved that for any prime $p\ge 5$,
\begin{align*}
\sum_{k=0}^{p-1}\frac{{4k\choose 2k}{2k\choose k}}{(2k+1)64^k}
\equiv (-1)^{\frac{p-1}{2}}-3p^2E_{p-3}\pmod{p^3},
\end{align*}
which was originally conjectured by Z.-W. Sun \cite{sunzw-scm-2011}.
Note that ${4k\choose 2k}/(2k+1)$ is the $2k$th Catalan number.

Mao and Z.-W. Sun \cite{ms-rj-2019} showed that for any prime $p\ge 5$,
\begin{align}
\sum_{k=0}^{(p-1)/2}\frac{{4k\choose 2k}{2k\choose k}}{(2k+1)64^k}\equiv (-1)^{\frac{p-1}{2}}2^{p-1}\pmod{p^2}.\label{new-7}
\end{align}
Z.-W. Sun \cite[Conjecture 11]{sunzw-2019} also conjectured an extension of \eqref{new-7} as follows.
\begin{conj}[Sun, 2019]
For any prime $p\ge 5$, we have
\begin{align}
\sum_{k=0}^{(p-1)/2}\frac{{4k\choose 2k}{2k\choose k}}{(2k+1)64^k}
\equiv (-1)^{\frac{p-1}{2}}\left(2^{p-1}-(2^{p-1}-1)^2\right)\pmod{p^3}.\label{a-3}
\end{align}
\end{conj}

The main purpose of the paper is to prove \eqref{a-3}. Our proof is based on hypergeometric series identities and symbolic summation method. Before proving \eqref{a-3}, we need the following two
key results.

\begin{thm}\label{t-1}
For any prime $p\ge 5$, we have
\begin{align}
&\sum_{k=0}^{(p-1)/2}\frac{{2k\choose k}}{(2k-1)^2 4^k}\equiv (-1)^{\frac{p-1}{2}}\left(2^{p-1}-q_p(2)\right)\pmod{p^2},\label{a-1}\\[10pt]
&\sum_{k=0}^{(p-1)/2}\frac{{2k\choose k}^2}{(2k-1)^3 16^k}\equiv 2-2q_p(2)-p\left(q_p(2)^2-4q_p(2)+3\right)\pmod{p^2},\label{a-2}
\end{align}
where $q_p(2)$ is the Fermat quotient $(2^{p-1}-1)/p$.
\end{thm}

\begin{thm}\label{t-2}
The supercongruence \eqref{a-3} is true.
\end{thm}

We shall prove Theorems \ref{t-1} and \ref{t-2} in Sections 2 and 3, respectively.

\section{Proof of Theorem \ref{t-1}}
\begin{lem}
For any integer $n\ge 2$, we have
\begin{align}
\sum_{k=0}^n\frac{(-n)_k (n-1)_k}{(1)_k\left(\frac{1}{2}\right)_k}=\frac{(-1)^{n-1}}{2n-1},\label{b-1}\\[10pt]
\sum_{k=0}^n\frac{(-n)_k (n-1)_k \left(-\frac{1}{2}\right)_k}{(1)_k^2\left(\frac{1}{2}\right)_k}
=\frac{4n(n-1)}{2n-1},\label{b-2}
\end{align}
where $(a)_0=1$ and $(a)_k=a(a+1)\cdots (a+k-1)$ for $k\ge 1$.
\end{lem}

{\noindent \it Proof.}
Recall Gauss's theorem \cite[(1.7.6), page 28]{slater-b-1966}:
\begin{align}
\pFq{2}{1}{a&b}{&c}{1}=\frac{\Gamma(c)\Gamma(c-a-b)}{\Gamma(c-a)\Gamma(c-b)},\label{new-3}
\end{align}
provided that $\Re(c-a-b)>0$. Letting $a=-n, b=n-1$ and $c=\frac{1}{2}$ in \eqref{new-3} gives
\begin{align*}
\pFq{2}{1}{-n&n-1}{&\frac{1}{2}}{1}=\frac{\Gamma\left(\frac{1}{2}\right)\Gamma\left(\frac{3}{2}\right)}
{\Gamma\left(\frac{1}{2}+n\right)\Gamma\left(\frac{3}{2}-n\right)}=\frac{(-1)^{n-1}}{2n-1},
\end{align*}
which is \eqref{b-1}.

Also, we have the following transformation formula of hypergeometric series \cite[(2.5.11), page 76]{slater-b-1966}:
\begin{align}
&\pFq{3}{2}{a&b&-n}{&e&f}{1}=\frac{(e-a)_n (f-a)_n}{(e)_n(f)_n}\notag\\[10pt]
&\times\pFq{3}{2}{1-s&a&-n}{&1+a-e-n&1+a-f-n}{1},\label{new-4}
\end{align}
where $s=e+f-a-b+n$. Letting $a=n-1,b=-\frac{1}{2},e=x$ and $f=\frac{3}{2}-x$ in \eqref{new-4} yields
\begin{align}
&\pFq{3}{2}{n-1&-\frac{1}{2}&-n}{&x&\frac{3}{2}-x}{1}=\frac{(x+1-n)_n \left(\frac{5}{2}-x-n\right)_n}{(x)_n\left(\frac{3}{2}-x\right)_n}\notag\\[10pt]
&\times\pFq{3}{2}{-2&n-1&-n}{&-x&x-\frac{3}{2}}{1}.\label{new-5}
\end{align}

Furthermore, we can evaluate the terminating hypergeometric series on the right-hand side of \eqref{new-5}:
\begin{align*}
&\pFq{3}{2}{-2&n-1&-n}{&-x&x-\frac{3}{2}}{1}\\[10pt]
&=\frac{4x^4-12x^3+(-8n^2+8n+11)x^2+(12n^2-12n-3)x+4n(n-1)(n^2-n-1)}{x(x-1)(2x-1)(2x-3)}.
\end{align*}
It follows that
\begin{align}
&\pFq{3}{2}{n-1&-\frac{1}{2}&-n}{&x&\frac{3}{2}-x}{1}=\frac{(x+1-n)_n \left(\frac{5}{2}-x-n\right)_n}{(x)_n\left(\frac{3}{2}-x\right)_n}\notag\\[10pt]
&\times\frac{4x^4-12x^3+(-8n^2+8n+11)x^2+(12n^2-12n-3)x+4n(n-1)(n^2-n-1)}{x(x-1)(2x-1)(2x-3)}.\label{new-6}
\end{align}
Letting $x\to 1$ on both sides of \eqref{new-6} and noting that
\begin{align*}
&\lim_{x\to 1}\frac{4x^4-12x^3+(-8n^2+8n+11)x^2+(12n^2-12n-3)x+4n(n-1)(n^2-n-1)}{x(2x-1)(2x-3)}\\
&=-4n^2(n-1)^2,
\end{align*}
and
\begin{align*}
\lim_{x\to 1}\frac{(x+1-n)_n \left(\frac{5}{2}-x-n\right)_n}{(x-1)(x)_n\left(\frac{3}{2}-x\right)_n}
=-\frac{1}{n(n-1)(2n-1)},
\end{align*}
we arrive at
\begin{align*}
\pFq{3}{2}{n-1&-\frac{1}{2}&-n}{&1&\frac{1}{2}}{1}=\frac{4n(n-1)}{2n-1},
\end{align*}
which proves \eqref{b-2}.
\qed

{\noindent \it Proof of \eqref{a-1}.}
We can rewrite \eqref{b-1} as
\begin{align}
\sum_{k=0}^{n-1}\frac{(-n)_k (n-1)_k}{(1)_k\left(\frac{1}{2}\right)_k}&=\frac{(-1)^{n-1}}{2n-1}
-\frac{(-n)_n (n-1)_n}{(1)_n\left(\frac{1}{2}\right)_n}\notag\\[10pt]
&=\frac{(-1)^{n-1}}{2n-1}\left(1+4^{n-1}(2n-2)\right)\notag\\[10pt]
&=(-1)^{n-1}\left(2^{2n-2}-\frac{2^{2n-2}-1}{2n-1}\right).\label{b-3}
\end{align}
Letting $n=\frac{p+1}{2}$ in \eqref{b-3} gives
\begin{align*}
\sum_{k=0}^{(p-1)/2}\frac{\left(\frac{-1-p}{2}\right)_k \left(\frac{-1+p}{2}\right)_k}
{(1)_k \left(\frac{1}{2}\right)_k}=(-1)^{\frac{p-1}{2}}\left(2^{p-1}-q_p(2)\right).
\end{align*}
Since for $0\le k \le \frac{p-1}{2}$,
\begin{align}
\left(\frac{-1-p}{2}\right)_k \left(\frac{-1+p}{2}\right)_k\equiv \left(-\frac{1}{2}\right)_k^2\pmod{p^2},\label{b-4}
\end{align}
we have
\begin{align}
\sum_{k=0}^{(p-1)/2}\frac{ \left(-\frac{1}{2}\right)_k^2}
{(1)_k \left(\frac{1}{2}\right)_k}\equiv (-1)^{\frac{p-1}{2}}\left(2^{p-1}-q_p(2)\right)\pmod{p^2}.\label{b-5}
\end{align}

Note that
\begin{align}
&\frac{\left(\frac{1}{2}\right)_k}{(1)_k}=\frac{{2k\choose k}}{4^k},\label{b-6}\\[10pt]
&\frac{\left(-\frac{1}{2}\right)_k}{\left(\frac{1}{2}\right)_k}=\frac{1}{1-2k}.\label{b-7}
\end{align}
Then the proof of \eqref{a-1} follows from \eqref{b-5}--\eqref{b-7}.
\qed

{\noindent \it Proof of \eqref{a-2}.}
We can rewrite \eqref{b-2} as
\begin{align}
\sum_{k=0}^{n-1}\frac{(-n)_k (n-1)_k \left(-\frac{1}{2}\right)_k}{(1)_k^2\left(\frac{1}{2}\right)_k}
&=\frac{4n(n-1)}{2n-1}-\frac{(-n)_n (n-1)_n \left(-\frac{1}{2}\right)_n}{(1)_n^2\left(\frac{1}{2}\right)_n}\notag\\[10pt]
&=\frac{1}{2n-1}\left(4n(n-1)+(-1)^n{2n-2\choose n}\right).\label{b-8}
\end{align}
Letting $n=\frac{p+1}{2}$ in \eqref{b-8} and using \eqref{b-4}, we obtain
\begin{align*}
\sum_{k=0}^{(p-1)/2}\frac{\left(-\frac{1}{2}\right)_k^3}{(1)_k^2\left(\frac{1}{2}\right)_k}
\equiv \frac{1}{p}\left(p^2-1+(-1)^{\frac{p+1}{2}}{p-1\choose \frac{p+1}{2}}\right)\pmod{p^2}.
\end{align*}

For $0\le k\le p-1$, we have
\begin{align}
{p-1\choose k}&\equiv (-1)^k\left(1-p\sum_{i=1}^k\frac{1}{i}+p^2\sum_{1\le i<j\le k}\frac{1}{ij}\right)\pmod{p^3}\notag\\[10pt]
&=(-1)^k\left(1-pH_k+\frac{p^2}{2}\left(H_k^2-H_k^{(2)}\right)\right),\label{new-1}
\end{align}
where
\begin{align*}
H_k^{(r)}=\sum_{j=1}^k\frac{1}{j^r},
\end{align*}
with the convention that $H_k=H_k^{(1)}$. It follows that
\begin{align}
\sum_{k=0}^{(p-1)/2}\frac{\left(-\frac{1}{2}\right)_k^3}{(1)_k^2\left(\frac{1}{2}\right)_k}
&\equiv \frac{p}{2}\left(H_{\frac{p+1}{2}}^2-H_{\frac{p+1}{2}}^{(2)}+2\right)-H_{\frac{p+1}{2}}\notag\\
&=\frac{p}{2}\left(H_{\frac{p-1}{2}}^2+4H_{\frac{p-1}{2}}-H_{\frac{p-1}{2}}^{(2)}+6\right)-H_{\frac{p-1}{2}}-2
\pmod{p^2}.\label{b-9}
\end{align}

By \cite[(41)]{lehmer-am-1938} and \cite[Lemma 2.4]{sunzw-scm-2011}, we have
\begin{align}
H_{\frac{p-1}{2}}\equiv -2q_p(2)+pq_p(2)^2\pmod{p^2},\label{b-10}
\end{align}
and
\begin{align}
H_{\frac{p-1}{2}}^{(2)}\equiv 0\pmod{p}.\label{b-11}
\end{align}
Substituting \eqref{b-10} and \eqref{b-11} into \eqref{b-9} gives
\begin{align}
\sum_{k=0}^{(p-1)/2}\frac{\left(-\frac{1}{2}\right)_k^3}{(1)_k^2\left(\frac{1}{2}\right)_k}\label{b-12}
\equiv 2q_p(2)-2+p\left(q_p(2)^2-4q_p(2)+3\right)\pmod{p^2}.
\end{align}

Finally, applying \eqref{b-6} and \eqref{b-7} to the left-hand side of \eqref{b-12}, we reach
\begin{align*}
\sum_{k=0}^{(p-1)/2}\frac{{2k\choose k}^2}{(2k-1)^3 16^k}\equiv
2-2q_p(2)-p\left(q_p(2)^2-4q_p(2)+3\right)\pmod{p^2},
\end{align*}
as desired.
\qed

\section{Proof of Theorem \ref{t-2}}
\begin{lem}
For any non-negative integer $n$, we have
\begin{align}
&\sum_{k=0}^n\frac{\left(-n\right)_k\left(n+1\right)_k\left(\frac{1}{4}\right)_k \left(\frac{3}{4}\right)_k}
{(1)_k^2\left(\frac{1}{2}\right)_k\left(\frac{3}{2}\right)_k}=\frac{{2n\choose n}}{4^n},\label{c-1}
\end{align}
and
\begin{align}
&\sum_{k=0}^n\frac{\left(-n\right)_k\left(n+1\right)_k\left(\frac{1}{4}\right)_k \left(\frac{3}{4}\right)_k}
{(1)_k^2\left(\frac{1}{2}\right)_k\left(\frac{3}{2}\right)_k}\sum_{j=1}^k\frac{1}{(2j-1)^2}\notag\\[10pt]
&=-\frac{{2n\choose n}}{4^n}\left(3+\sum_{k=1}^n\frac{1}{(2k-1)^2}\right)+\frac{2}{2n+1}\sum_{k=0}^n
\frac{{2k\choose k}}{(2k-1)^2 4^k}\notag\\[10pt]
&-\frac{4^n}{(2n+1){2n\choose n}}\sum_{k=0}^n
\frac{{2k\choose k}^2}{(2k-1)^3 16^k}.\label{c-2}
\end{align}
\end{lem}
{\noindent \it Proof.} Recall that (see \cite[(2.4.2.2), page 65]{slater-b-1966})
\begin{align}
\pFq{4}{3}{d&1+f-g&\frac{f}{2}&\frac{f+1}{2}}{&1+f&\frac{1+f+d-g}{2}&1+\frac{f+d-g}{2}}{1}
=\frac{\Gamma(g-f)\Gamma(g-d)}{\Gamma(g)\Gamma(g-f-d)}.\label{new-2}
\end{align}
Letting $d=-n,f=\frac{1}{2}$ and $g=-n+\frac{1}{2}$ in \eqref{new-2}, we obtain
\begin{align*}
\pFq{4}{3}{-n&n+1&\frac{1}{4}&\frac{3}{4}}{&1&\frac{1}{2}&\frac{3}{2}}{1}
=\frac{\Gamma(-n)\Gamma\left(\frac{1}{2}\right)}{\Gamma\left(-n+\frac{1}{2}\right)\Gamma(0)}
=\frac{{2n\choose n}}{4^n},
\end{align*}
which is \eqref{c-1}.

On the other hand, \eqref{c-2} can be discovered and proved by symbolic summation
package {\tt Sigma} due to Schneider \cite{schneider-slc-2007}. One can refer to \cite{liu-jsc-2019} for the same approach to finding and proving identities of this type.
\qed

{\noindent \it Proof of \eqref{a-3}.}
Recall that (see \cite[(4.4)]{liu-jmaa-2019})
\begin{align}
\left(\frac{1+p}{2}\right)_k\left(\frac{1-p}{2}\right)_k
\equiv \left(\frac{1}{2}\right)_k^2\left(1-p^2\sum_{j=1}^k\frac{1}{(2j-1)^2}\right)
\pmod{p^4}.\label{c-3}
\end{align}
Letting $n=\frac{p-1}{2}$ in \eqref{c-1} and using \eqref{c-3}, we obtain
\begin{align}
\sum_{k=0}^{(p-1)/2}\frac{\left(\frac{1}{2}\right)_k\left(\frac{1}{4}\right)_k \left(\frac{3}{4}\right)_k}{(1)_k^2\left(\frac{3}{2}\right)_k}
&\equiv \frac{1}{2^{p-1}}{p-1\choose \frac{p-1}{2}}\notag\\
&+p^2\sum_{k=0}^{(p-1)/2}\frac{\left(\frac{1}{2}\right)_k\left(\frac{1}{4}\right)_k \left(\frac{3}{4}\right)_k}{(1)_k^2\left(\frac{3}{2}\right)_k}\sum_{j=1}^k\frac{1}{(2j-1)^2}
\pmod{p^4},\label{c-4}
\end{align}
where we have utilized the fact $\left(\frac{1}{4}\right)_k \left(\frac{3}{4}\right)_k/\left(\frac{3}{2}\right)_k \in \mathbb{Z}_p$ for $0\le k\le \frac{p-1}{2}$.

From \eqref{c-3}, we deduce that
\begin{align}
\left(\frac{1+p}{2}\right)_k\left(\frac{1-p}{2}\right)_k
\equiv \left(\frac{1}{2}\right)_k^2\pmod{p^2}.\label{c-5}
\end{align}
Letting $n=\frac{p-1}{2}$ in \eqref{c-2} and using \eqref{c-5} gives
\begin{align}
&\sum_{k=0}^{(p-1)/2}\frac{\left(\frac{1}{2}\right)_k\left(\frac{1}{4}\right)_k \left(\frac{3}{4}\right)_k}{(1)_k^2\left(\frac{3}{2}\right)_k}\sum_{j=1}^k\frac{1}{(2j-1)^2}\notag\\[10pt]
&\equiv -\frac{1}{2^{p-1}}{p-1\choose \frac{p-1}{2}}\left(3+\sum_{k=1}^{(p-1)/2}\frac{1}{(2k-1)^2}\right)+\frac{2}{p}\sum_{k=0}^{(p-1)/2}
\frac{{2k\choose k}}{(2k-1)^2 4^k}\notag\notag\\[10pt]
&-\frac{2^{p-1}}{p{p-1\choose \frac{p-1}{2}}}\sum_{k=0}^{(p-1)/2}
\frac{{2k\choose k}^2}{(2k-1)^3 16^k}\pmod{p^2}.\label{c-6}
\end{align}
Substituting \eqref{c-6} into \eqref{c-4} yields
\begin{align}
\sum_{k=0}^{(p-1)/2}\frac{\left(\frac{1}{2}\right)_k\left(\frac{1}{4}\right)_k \left(\frac{3}{4}\right)_k}{(1)_k^2\left(\frac{3}{2}\right)_k}
&\equiv \frac{1}{2^{p-1}}{p-1\choose \frac{p-1}{2}}
-\frac{p^2}{2^{p-1}}{p-1\choose \frac{p-1}{2}}\left(3+\sum_{k=1}^{(p-1)/2}\frac{1}{(2k-1)^2}\right)\notag\\[10pt]
&+2p\sum_{k=0}^{(p-1)/2}\frac{{2k\choose k}}{(2k-1)^2 4^k}
-\frac{2^{p-1}p}{{p-1\choose \frac{p-1}{2}}}\sum_{k=0}^{(p-1)/2}
\frac{{2k\choose k}^2}{(2k-1)^3 16^k}\pmod{p^4}.\label{c-7}
\end{align}

Furthermore, by \eqref{new-1}, \eqref{b-10} and \eqref{b-11} we have
\begin{align}
{p-1\choose \frac{p-1}{2}}&\equiv (-1)^{\frac{p-1}{2}}\left(1-pH_{\frac{p-1}{2}}+\frac{p^2}{2}\left(H_{\frac{p-1}{2}}^2-H_{\frac{p-1}{2}}^{(2)}\right)\right)\notag\\
&\equiv (-1)^{\frac{p-1}{2}}\left(1+2pq_p(2)+p^2q_p(2)^2\right)\pmod{p^3}.\label{c-8}
\end{align}
By \eqref{b-11} and the Wolstenholme's theorem \cite[page 114]{hw-b-2008}, we have
\begin{align}
\sum_{k=1}^{(p-1)/2}\frac{1}{(2k-1)^2}=H_{p-1}^{(2)}-\frac{1}{4}H_{\frac{p-1}{2}}^{(2)}
\equiv 0\pmod{p}.\label{c-9}
\end{align}

Setting $2^{p-1}=a$ and $q_p(2)=(a-1)/p$, and then
substituting \eqref{a-1}, \eqref{a-2}, \eqref{c-8} and \eqref{c-9} into \eqref{c-7}, we arrive at
\begin{align*}
\sum_{k=0}^{(p-1)/2}\frac{\left(\frac{1}{2}\right)_k\left(\frac{1}{4}\right)_k \left(\frac{3}{4}\right)_k}{(1)_k^2\left(\frac{3}{2}\right)_k}
&\equiv (-1)^{\frac{p-1}{2}}\left(\frac{a^3-2a^2+4a-2}{2a-1}+\frac{3(a-1)^2p^2}{a(2a-1)}\right)\pmod{p^3}\\
&=(-1)^{\frac{p-1}{2}}\left(a-(a-1)^2+\frac{3(a-1)^3}{2a-1}+\frac{3(a-1)^2p^2}{a(2a-1)}\right).
\end{align*}
By the Fermat's little theorem, we have $a-1\equiv 0\pmod{p}$, and so
\begin{align}
\sum_{k=0}^{(p-1)/2}\frac{\left(\frac{1}{2}\right)_k\left(\frac{1}{4}\right)_k \left(\frac{3}{4}\right)_k}{(1)_k^2\left(\frac{3}{2}\right)_k}\equiv
(-1)^{\frac{p-1}{2}}\left(a-(a-1)^2\right)\pmod{p^3}.\label{c-10}
\end{align}

Note that
\begin{align}
&\frac{\left(\frac{1}{4}\right)_k\left(\frac{3}{4}\right)_k}{(1)_k^2}=\frac{{4k\choose 2k}{2k\choose k}}{64^k},\label{c-11}\\[10pt]
&\frac{\left(\frac{1}{2}\right)_k}{\left(\frac{3}{2}\right)_k}=\frac{1}{2k+1}.\label{c-12}
\end{align}
Then the proof of \eqref{a-3} follows from \eqref{c-10}--\eqref{c-12}.
\qed

\vskip 5mm \noindent{\bf Acknowledgments.}
This work was supported by the National Natural Science Foundation of China (grant 11801417).

\end{document}